\documentclass{amsart}
\chardef\bslash=`\\

\usepackage{amssymb,amsmath,amsfonts,graphicx,amsthm,amscd,stmaryrd,mathrsfs,color,mathdots}

\usepackage[all,cmtip,poly]{xy}
\newcommand*{\Scale}[2][4]{\scalebox{#1}{$#2$}}

\newtheorem{theorem}[subsection]{Theorem}

\newtheorem{lemma}[subsection]{Lemma}

\newtheorem{proposition}[subsection]{Proposition}

\newtheorem{defn}[subsection]{Definition}
\theoremstyle{remark}
\newtheorem{remark}[subsection]{Remark}

\numberwithin{equation}{section}

\newif\iffinalrun
\iffinalrun
\else
 \fi

\iffinalrun
  \newcommand{\need}[1]{}
  \newcommand{\mar}[1]{}
\else
  \newcommand{\need}[1]{{\tiny *** #1}}
  \newcommand{\mar}[1]{\marginpar{\raggedright\tiny  #1}}\fi

\renewcommand\mathbb{\mathbf}

\newcommand{\Lie}{{\operatorname{Lie}\,}}

\newcommand{\frg}{{\mathfrak{g}}}
\newcommand{\frs}{{\mathfrak{s}}}
\newcommand{\frl}{{\mathfrak{l}}}

\renewcommand{\ell}{l}

\newcommand{\diag}{\operatorname{diag}}

\newcommand{\bC}{\ensuremath{\mathbf{C}}}

\newcommand{\bQ}{\ensuremath{\mathbf{Q}}}

\newcommand{\bR}{\ensuremath{\mathbf{R}}}
\newcommand{\bbC}{\ensuremath{\mathbf{C}}}

\newcommand{\bZ}{\ensuremath{\mathbf{Z}}}

\newcommand{\bbZ}{\ensuremath{\mathbf{Z}}}
\newcommand{\bbQ}{\ensuremath{\mathbf{Q}}}
\newcommand{\bbR}{\ensuremath{\mathbf{R}}}

\newcommand{\cF}{{\mathcal F}}

\newcommand{\cR}{{\mathcal R}}

\newcommand{\frh}{{\mathfrak h}}

\newcommand{\frc}{\mathfrak{c}}

\DeclareMathOperator{\Aut}{Aut}
\DeclareMathOperator{\Ad}{Ad}

\newcommand{\GL}{\mathrm{GL}}

\DeclareMathOperator{\SL}{SL}
\DeclareMathOperator{\Spec}{Spec}

\DeclareMathOperator{\Sel}{Sel}

\newcommand{\Res}{\operatorname{Res}}
\newcommand{\Pic}{\operatorname{Pic}}

\newcommand{\SO}{\operatorname{SO}}

\newcommand{\doubleslash}{/\kern-0.2em{/}}

\begin{document}

\author{Jack A. Thorne}
\title[Reduction theory]{Reduction theory for stably graded Lie algebras}
\begin{abstract} We define a reduction covariant for the representations \`a la Vinberg associated to stably graded Lie algebras. We then give an analogue of the LLL algorithm for the odd split special orthogonal group and show how this can be combined with our theory to effectively reduce the coefficients of vectors in a representation connected to 2-descent for odd hyperelliptic curves. 
\end{abstract}
\maketitle
\setcounter{tocdepth}{1}
\tableofcontents
\section{Introduction}

Let $\frh$ be a semisimple Lie algebra over a field $k$ of characteristic $0$, and let $m \geq 1$. By a $\bbZ / m \bbZ$-grading, we mean a direct sum decomposition
\[ \frh = \oplus_{i \in \bbZ / m \bbZ} \frh_i, \]
where for all $i, j$ we have $[\frh_i, \frh_j] \subset \frh_{i + j}$. To such a grading, we can associate a pair $(G, V)$, where $G$ is the connected component of the identity in the group of automorphisms of $\frh$ which preserve the grading, and $V = \frh_1$ is a $k$-vector space which is naturally a representation of the reductive group $G$. This representation is coregular, in the sense that  the ring $k[V]^G$ of invariant polynomials is freely generated and the geometric quotient $B = \Spec k[V]^G$ is isomorphic to affine space \cite{Vin76}.

We say that the grading is stable if $V$ contains $G$-orbits which are stable in the sense of geometric invariant theory (i.e.\ closed, with finite stabilizers in $G$). In this case the $G$-invariant locus $V^s \subset V$ of stable vectors is Zariski open and non-empty, and it is of interest to study the set of $G(k)$-orbits in a given geometric stable orbit. The stable gradings have been classified \cite{Ree12} (at least when $k$ is algebraically closed) and experience shows that it is frequently the case that there is a family of projective curves $C \to B^s$ and for each $b \in B^s(k)$ an injection
\[ \Pic^0 C_b(k) / m \Pic^0 C_b(k) \hookrightarrow G(k) \backslash V_b(k). \]
When $k = \bbQ$, this can often be extended to an injection (with source the $m$-Selmer group of the Jacobian variety $ \Pic^0 C_b$)
\[ \Sel_m \Pic^0 C_b \hookrightarrow G(\bbQ) \backslash V_b(\bbQ), \]
with image contained in the image of the map
\[ G(\bbZ) \backslash V_b(\bbZ) \to G(\bbQ) \backslash V_b(\bbQ) \]
(for a suitable choice of integral structures). For example, this is the case for  stable $\bbZ / m \bbZ$-gradings of the semisimple Lie algebras of types $A_2$, $D_4$, $E_6$ and $E_8$ for $m = 2, 3, 4$ and $5$, respectively, in which case $C$ can be taken to be the Weierstrass family of elliptic curves (see \cite{Cre10, Fis13}), and for the stable $\bbZ / 2 \bbZ$-gradings of the semisimple Lie algebras of types $A_n$, $D_n$, and $E_n$ (see e.g. \cite{Bha13, Sha19, Lag22}). 

It is therefore of particular interest to be able to understand the sets $G(\bbZ) \backslash V_b(\bbZ)$ of integral orbits. This is the goal of reduction theory, which can be understood in this context to have two steps.

In the first step, one defines a \emph{reduction covariant}, i.e.\ a $G(\bbR)$-equivariant map
\[ \cR : V^s(\bbR) \to X_G, \]
where $X_G$ is the symmetric space of $G$ (i.e. the homogeneous space $G(\bbR) / K$, where $K$ is a maximal compact mod centre subgroup of $G(\bbR)$). In the second step, one introduces a notion of reducedness for elements of the symmetric space $X_G$, and says that an element $T \in V^s(\bbR)$ is reduced if its reduction covariant $\cR(T)$ satisfies this condition. For example, suppose given a fundamental set $\cF \subset X_G$ for the action of $G(\bbZ)$, i.e.\ a subset of $\cF$ such that $G(\bbZ) \cdot \cF = X_G$ and $\cF$ intersects each orbit of $G(\bbZ)$ finitely many times. Then we could say that $T$ is $\cF$-reduced if $\cR(T) \in \cF$. An immediate consequence is that any $T \in V^s(\bbR)$ is $G(\bbZ)$-conjugate to an $\cF$-reduced element, and that if $b \in B(\bbZ) \cap B^s(\bbR)$ then $V_b(\bbZ)$ contains only finitely many $\cF$-reduced elements. 

In this paper we consider both the first and second steps described in the previous paragraph. Reduction covariants have been defined in the literature in some isolated cases for representations arising from stably graded Lie algebras (including, in the papers \cite{Cre10, Fis13},  for the representations $V$ associated to $m$-descent on elliptic curves for $m = 2, 3, 4, 5$). The first main result (Theorem \ref{thm_special_cartan_involution}) of this paper is the definition of a reduction covariant $ \cR : V^s(\bbR) \to X_G$ for any stably $\bbZ / m \bbZ$-graded Lie algebra. To define it, we think of $X_G$ as the space of Cartan involutions of the real semisimple group $G_\bbR$.  It turns out that for any  Cartan subspace $\frc \subset V_\bbR$, there is a unique Cartan involution of $\frh_\bbR$ which respects the grading and leaves $\frc$ invariant. A stable vector $T$ is contained in a unique Cartan subspace, and we define $\cR(T)$ to be the restriction of the Cartan involution  associated to this subspace to $G_\bbR$. 

Regarding the second step, a good starting point is the reduction theory of Borel and Harish-Chandra \cite{Bor62, Bor19}. They defined Siegel sets $\mathfrak{S}_t \subset G(\bbR)$, which are rather explicit fundamental sets in the sense defined above, and which can serve to define a notion of reducedness. However, what is lacking in general is an explicit algorithm to transport elements into the Siegel set. One case where the existing theory works well is when $G$ is isogenous to $\SL_n$ (or more generally, a product of such groups). In this case $X_G$ may be identified with the set of lattices in $\bbR^n$ of covolume 1 and lattice reduction algorithms (such as the LLL algorithm described in \cite{Len82}) apply without change. The second main result of this paper, which appears in \S \ref{sec_reduction_algorithm}, is an algorithm to transport elements of $X_G$ into a Siegel set in the case that $G$ is a split special orthogonal group $\SO_{2g+1}$. This is the group that appears for the stable $\bbZ / 2 \bbZ$-grading of $\frs\frl_{2g+1}$, used in \cite{Bha13} to study the 2-Selmer groups of Jacobians of odd hyperelliptic curves. This algorithm is based on the LLL algorithm. Roughly speaking, we observe that the two main steps of the LLL algorithm correspond either to acting by ``integral unipotent transformations'' or ``simple reflections in $S_n$''. We simply replace $S_n$ by the Weyl group of $\SO_{2g+1}$. We consider the analogue of this, where $\SO_{2g+1}$ is replaced by any split semisimple group, in another paper \cite{Rom24}.

\subsection{Outline} In \S \ref{sec_covariant}, we define the reduction covariant $\cR$ using the theory of the Cartan involution. In \S \ref{sec_example_grading}, we compute the covariant explicitly in the case of the stable $\bbZ / 2 \bbZ$-grading of $\frs\frl_{2g+1}$. In \S \ref{sec_reduction_algorithm}, we give our algorithm to move a given element of $X_{\SO_{2g+1}}$ into a Siegel set. Finally, in \S \ref{sec_example_computation}, we show how these ideas work together in an explicit example. 

\subsection{Notation}

In this paper, a reductive group over a field $k$ means a  smooth connected linear algebraic group over a field $k$ of trivial unipotent radical. In particular, by a reductive group $H$ over $\bbR$ we mean a connected algebraic group (although the associated set $H(\bbR)$ of real points might not be connected). Similarly, a semisimple group is a reductive group of trivial radical. 

We will use gothic letters to denote Lie algebras (so $\Lie H = \mathfrak{h}$) and subscripts to denote base extension (so if $H$ is an algebraic group over $\bbR$, then $(\Lie H) \otimes_\bR \bC = \mathfrak{h}_\bC$).  We use the notations $Z_H(\cdot)$ and $Z_\frh(\cdot)$ to denote group and Lie algebra centralizer, respectively.

We use superscript $0$ to denote the connected component of the identity. Thus if $H$ is a linear algebraic group over $\bbR$ then there is an inclusion $H(\bR)^0 \subset H^0(\bR)$ which is not in general an equality.
\section{The reduction covariant}\label{sec_covariant}

We will summarise the necessary properties of Cartan involutions of real reductive groups before proving the existence of the reduction covariant described in the introduction.

\subsection{Background on Cartan involutions}

\begin{defn}
Let $\frh$ be a semisimple Lie algebra over $\bbR$. A Cartan involution of $\frh$ is a Lie algebra involution $\theta : \frh \to \frh$ satisfying the following equivalent conditions:
\begin{enumerate}
\item The symmetric bilinear form $B_\theta(X, Y) = -B(X, \theta Y) : \frh \times \frh \to \bbR$ (where $B$ is the Killing form) is positive definite.
\item The $\bbR$-vector subspace $\{ X \in \frh_\bC \mid \theta(X) = \overline{X} \}$ is a real form of $\frh_\bC$ which is compact, in the sense that its Killing form is negative definite.
\end{enumerate}
\end{defn}
It is convenient to define Cartan involutions for disconnected groups. We do this following \cite{Ada18} (see also \cite{Mos55}).
\begin{defn}
Let $H$ be a linear algebraic group over $\bbR$ such that $H^0$ is reductive. A Cartan involution of $H$ is an involution $\theta : H \to H$ such that
\[ K_\theta = \{ h \in H(\bbC) \mid \theta(h) = \overline{h} \} \]
is a compact subgroup of $H(\bbC)$ that meets every connected component of $H(\bbC)$. 
\end{defn}
If $\theta$ is a Cartan involution, then $H(\bbR)^\theta = K_\theta \cap H(\bbR)$ is a maximal compact subgroup of $H(\bbR)$ (\cite[\S 4]{Ada18}). 

Given an automorphism $\theta$ of an algebraic group $H$, we will also write $\theta$ for the induced automorphism of its Lie algebra $\frh$. The following proposition is well-known. 
\begin{proposition}
Let $H$ be a semisimple group over $\bbR$, and let $\theta$ be an involution of $H$. Then $\theta : H \to H$ is a Cartan involution if and only if $\theta : \frh \to \frh$ is a Cartan involution. 
\end{proposition}

\begin{proposition}\label{prop_facts_about_Cartan_involutions} Let $H$ be a reductive group over $\bbR$.
\begin{enumerate} 
\item Cartan involutions of $H$ exist. If $\theta, \theta'$ are Cartan involutions of $H$, then there exists $h \in H(\bbR)^0$ such that $\theta' = \Ad(h) \circ \theta \circ \Ad(h^{-1})$. 
\item Suppose $M \subset H$ is a closed subgroup and $\theta$ is a Cartan involution of $H$ such that $\theta(M) = M$. Then $M^0$ is reductive and $\theta|_M$ is a Cartan involution of $M$.
\item Suppose $M \subset H$ is a closed subgroup with  $M^0$ reductive and that $\theta_M$ is a Cartan involution of $M$. Then there exists a Cartan involution $\theta_H$ of $H$ extending $\theta_M$.
\end{enumerate}
\end{proposition}
\begin{proof}
See \cite[Theorem 3.12]{Ada18}.
\end{proof} 
If $H$ is a reductive group over $\bbR$, then we write $X_H$ for the set of Cartan involutions of $H$. Proposition \ref{prop_facts_about_Cartan_involutions}(1) shows that this is a homogeneous space for $H(\bbR)$. 

\subsection{Application to graded Lie algebras}
\begin{theorem}\label{thm_special_cartan_involution}
Let $\frh$ be a semisimple Lie algebra over $\bbR$. Suppose given a stable $\bbZ / m \bbZ$-grading $\frh = \oplus_{i \in \bbZ / m \bbZ} \frh_i$, and let $(G, V)$ be the associated pair. Then:
\begin{enumerate}
\item For any stable vector $T \in V^s(\bbR)$, there exists a unique Cartan involution $\theta_T$ of $H$ such that  $[T, \theta_T(T)] = 0$ and $\theta_T(\frh_i) \subset \frh_{-i}$ for each $i \in \bbZ / m \bbZ$.
\item The map $\cR = \cR_\frh : V^s(\bbR)\to X_G$ defined by $T \mapsto \theta_T|_G$ is $G(\bbR)$-equivariant.
\end{enumerate}
\end{theorem}
\begin{proof}
 Let $H$ denote the adjoint group over $\bbR$ with Lie algebra $\frh$; then $\Aut(H) = \Aut(\frh)$, and the embedding $\Ad : H \to \Aut(H)$ allows us to identify $H$ with the connected component of the identity of the linear algebraic group $\Aut(H)$.  Since $T$ is stable, it is regular semisimple when viewed as an element of $\frh$ (\cite[Lemma 13]{Ree12}). Let  $\frc = Z_\frh(T)$, a Cartan subalgebra of $\frh$ since $T$ is regular semisimple, and let $C = Z_H(T)$ denote the maximal torus of $H$ with Lie algebra $\frc$.

Let us first establish uniqueness. If a Cartan involution $\theta_T$ with the given properties exists, then the identity $[T, \theta_T(T)] = 0$ shows that $\theta_T$ normalises $\frc$ and $C$. The restriction of $\theta_T$ to $C$ must be a Cartan involution of $C$. The torus $C$ has a unique Cartan involution $\theta_C$, so if $\theta_T'$ is another Cartan involution with the given properties, then we can write $\theta_T' = t \theta_T$ for some $t \in C(\bbR)$. 

Giving an $\bbZ / m \bbZ$-grading is equivalent to giving a homomorphism $\mu_m \to \Aut(\frh)$. Let $\sigma : \frh_\bC \to \frh_\bC$ be the image of $e^{2 \pi i / m} \in \mu_m(\bC)$. The inclusion $\theta_T(\frh_i) \subset \frh_{-i}$ may be equivalently expressed as $\sigma \theta_T = \theta_T \overline{\sigma}$. The identities $\sigma \theta_T = \theta_T \overline{\sigma}$, $\sigma t \theta_T = t \theta_T \overline{\sigma}$ together imply that $t \in C^\sigma(\bbR)$. This group is finite (because the grading is stable), hence compact, which implies that $t, \theta_T$ commute, and so $t^2 = 1$. We now consider the positive definite symmetric bilinear forms $B_{\theta_T}$, $B_{t \theta_T}$ on $\frh$. Computing
\[ B_{t \theta_T}(X, Y) = -B(X, t \theta_T Y) = -B(t^{-1}X, \theta_T Y) = B_{\theta_T}(t^{-1} X, Y) = B_{\theta_T}(X,  t Y) \]
shows that (with respect to the inner product determined by $B_{\theta_T}$) $t$ is an orthogonal endomorphism of $\frh$ with the property that $B_{\theta_T}(t X, X) > 0$ for every non-zero $X \in \frh$. Since $t^2 = 1$, this is only possible if $t = 1$, hence $\theta_T = \theta_T'$.

We now establish existence of $\theta_T$. Let $M = C \rtimes \mu_m \subset \Aut(H)$, and let $\theta_M$ be the Cartan involution of $M$ given by the formula $\theta_M(c \rtimes \sigma) = \theta_C(c) \rtimes \overline{\sigma}$. To show that this is a group homomorphism, we must show that the identity $\sigma \theta_C = \theta_C \overline{\sigma}$ holds in $\Aut(C_\bC)$.  Let $K_C \subset C(\bbC)$ be the unique maximal compact subgroup. If $x \in K_C$, then $\theta_C(x) = \overline{x}$ and $\sigma(x) \in K_C$. For such $x$, we have $\sigma \theta_C(x) = \sigma(\overline{x})$ and $\theta_C \overline{\sigma}(x) = \overline{\overline{\sigma}(x)} = \sigma(\overline{x}) = \sigma \theta_C(x)$. Since $K_C$ is Zariski dense in $C_\bC$, this shows that in fact $\sigma \theta_C = \theta_C \overline{\sigma}$.

By Proposition \ref{prop_facts_about_Cartan_involutions}, we can extend $\theta_M$ to a Cartan involution $\theta$ of $\Aut(H)$. We set $\theta_T = \theta|_H$. We claim that $\theta_T$ has the desired properties. It normalises $C$ by construction, so satisfies $[T, \theta_T(T)] = 0$. We must show that $\sigma \theta_T = \theta_T \overline{\sigma}$ in $\Aut(H)(\bC)$. To this end, we note that if $h \in H(\bC) \subset \Aut(H)(\bC)$, and $\alpha \in \Aut(H)(\bC)$, then $\alpha h \alpha^{-1} = \alpha(h)$. (More transparently, we have $\alpha \circ \Ad(h) \circ \alpha^{-1} = \Ad(\alpha(h))$.) Since $\theta : \Aut(H) \to \Aut(H)$ is a group homomorphism, we have
\[ \theta(\alpha h \alpha^{-1}) = \theta(\alpha) \theta(h) \theta(\alpha)^{-1}. \]
The right-hand side equals $\theta(\alpha)(\theta_T(h))$, while the left-hand side equals
\[ \theta( \alpha(h) ) = \theta_T( \alpha(h) ). \]
Since $h$ is arbitrary, this gives the identity $\theta(\alpha) \circ \theta_T = \theta_T \circ \alpha$ (composition in $\Aut(H)(\bC)$). Taking $\alpha = \overline{\sigma}$ and using the identity $\theta(\overline{\sigma}) = \theta_M(\overline{\sigma}) = \sigma$, we find $\sigma \theta_T= \theta_T \overline{\sigma} $, as required.
This completes the proof of the first part of the Theorem. The second part follows from the properties which uniquely characterize $\theta_T$. 
\end{proof}
We now record some basic properties of the reduction covariant $\cR$.
\begin{proposition}\label{prop_base_extension}
\begin{enumerate} \item Let $\frh, \frh'$ be semisimple Lie algebras over $\bbR$. Suppose given stable $\bbZ / m \bbZ$-gradings 
\[ \frh = \oplus_{i \in \bbZ / m \bbZ} \frh_i, \text{ } \frh' = \oplus_{i \in \bbZ / m \bbZ} \frh'_i. \]
Let $\frh'' = \frh \oplus \frh'$ with its induced $\bbZ / m \bbZ$-grading. Then, with the obvious notation, we have $\cR_{\frh''} = \cR_{\frh} \times \cR_{\frh'}$. 
\item Let $\frh$ be a semisimple Lie algebra over $\bbR$. Suppose given a stable $\bbZ / m \bbZ$-grading $\frh = \oplus_{i \in \bbZ / m \bbZ} \frh_i$, and let $(G, V)$ be the associated pair. Let $\frh' = \Res_{\bC / \bR} (\frh_\bC)$, with its induced grading $\frh' = \oplus_{i  \in \bbZ / m \bbZ} \Res_{\bC / \bR} (\frh_{i, \bC})$. Then this grading is stable and its associated pair may be identified 
as
\[ (G', V') = (\Res_{\bC / \bR} G_\bC, \Res_{\bC / \bR} V_\bC). \]
 Moreover, there is a commutative diagram
\[ \xymatrix{ V^s(\bbR) \ar[d]_{\cR_\frh} \ar[r] & V^s(\bC) = (V')^s(\bR) \ar[d]^{\cR_{\frh'}} \\
X_G \ar[r] &   X_{G'} }. \]
\end{enumerate} 
\end{proposition}
\begin{proof}
The first part follows quickly from the definition of $\cR$. For the second, we note that the stability of the grading of $\frh'$ can be checked over $\bC$. There is an isomorphism $\frh'_\bC \cong \frh_\bC \oplus \frh_\bC$, so the grading of $\frh'$ is stable. It is easy to check that $(G', V')$ has the claimed form. 

All arrows in the commutative diagram are the natural ones, except for the bottom arrow, which we need to define. Let $\theta_G$ be a Cartan involution of $G$. We can then check directly from the definition that the involution $\theta'_G$ of $G'_\bC \cong G_\bC \times G_\bC$, given on complex points $(x_1, x_2) \in G(\bC) \times G(\bC)$ by the formula $\theta'_G(x_1, x_2) = (\theta_G(x_2), \theta_G(x_1))$, is defined over $\bR$ and is a Cartan involution. We send $\theta_G$ to $\theta'_G$. 

To show that the diagram is commutative, let $T \in V^s(\bbR)$ and let $\theta_T$ be the associated Cartan involution of $H$, the adjoint group of $\frh$. Let $T' \in V^s(\bC) = (V')^s(\bR)$ be the image of $T$ under the natural map. It follows from the definitions that $ \theta_T'$ (where $\theta_T'$ is given by the same formula as in the previous paragraph) satisfies the conditions characterizing $\theta_{T'}$. This completes the proof. 
\end{proof} 
One possible use of the first part of Proposition \ref{prop_base_extension} is the case where we start with a stably graded Lie algebra $\frh_0$ over $\bQ$, and take $\frh = (\Res_{K / \bQ} \frh_0)_\bR$, for a number field $K / \bQ$. We find that the reduction covariant of $V^s(\bR) = V_0^s(K \otimes_\bQ \bR)$ is the product (over the set of infinite places $v$ of $K$) of the reduction covariants for the stably graded Lie algebras $\Res_{K_v / \bR} \frh_{0, K_v}$. One possible use of the second part is in computing the reduction covariant of $V^s(\bR)$ explicitly. For example, if we start with the stable $\bbZ / 3 \bbZ$-grading of the exceptional Lie algebra $\frg_2$, then $V$ is (close to) the space of binary cubic forms. In \cite{Sto03}, several reduction covariants on the space of binary cubic forms are considered, but it is shown that there is a unique one, namely the Julia invariant, which is compatible with extension of scalars from $\bbR$ to $\bbC$ (see \cite[Proposition 3.4]{Sto03}). This characterisation may be used to relate the reduction covariant $\cR$ we construct in this case to the Julia invariant. 

\section{Example: The stable $\bbZ / 2 \bbZ$-grading of $\frs\frl_{2g+1}$}\label{sec_example_grading}

The stable $\bbZ / m \bbZ$-gradings of semisimple Lie algebras over $\bbC$ have been classified \cite{Ree12}. The first case is when $m = 2$; then each semisimple Lie algebra over $\bbC$ has a unique stable $\bbZ / 2 \bbZ$-grading, up to isomorphism. 

Suppose instead that $\frh$ is a semisimple Lie algebra over a field $k$ of characteristic $0$, and let $H$ be its associated adjoint group. In this case there need not be a unique $H(k)$-conjugacy class of stable $\bbZ / 2 \bbZ$-gradings over $k$, but if $H$ is split then there is a unique $H(k)$-conjugacy class of stable $\bbZ / 2 \bbZ$-gradings with the property that $\frh_1$ contains a regular nilpotent element (see \cite[Corollary 2.14]{Tho13}). When $k = \bbQ$, this distinguished class of stable gradings has been used to study the 2-Selmer groups of families of algebraic curves \cite{Bha13, Tho15, Sha19, Lag22}. In this section we describe this distinguished $\bbZ / 2 \bbZ$-grading in Dynkin type $A_{2g}$ and make explicit the construction of Theorem \ref{thm_special_cartan_involution}. 

Let $g \geq 1$ be an integer, and let $e_{-g}, \dots, e_{-1}, e_0, e_1, \dots, e_g$ denote the standard basis of $\bbR^{2g+1}$ (with a shift in the indices). Let $J$ be the Gram matrix of the symmetric bilinear form defined by $\langle e_i, e_{-j} \rangle = \delta_{ij}$ if $0 \leq i, j \leq g$, $\langle e_i, e_j \rangle = 0$ if $1 \leq i, j \leq g$. Thus
\[ J = \left( \begin{array}{ccc} & & 1 \\
& \iddots &  \\
1 \end{array}\right). \]
We define an involution of $\frh = \frs\frl_{2g+1}$ (equivalently, a $\bbZ / 2 \bbZ$-grading) by the formula $\sigma(T) = - J {}^t T J$. This is stable. The associated group is 
\[ G = \mathrm{SO}(J) = \{ x \in \SL_{2g+1} \mid {}^t x J x = J\}, \]
 and the associated representation is 
 \[ V = \{ T \in \frs\frl_{2g+1} \mid J {}^t T J = T \}. \]
 The open subscheme $V^s \subset V$ consists of those operators whose characteristic polynomial has no repeated roots.  
 
 The following proposition describes the map $\cR : V^s(\bR) \to X_G$ in terms of linear algebra.
\begin{proposition}\label{prop_explicit_reduction_theory}
\begin{enumerate}
\item We can identify $X_G$ with the set of inner products $H$ on $\bbR^{2g+1}$ which are compatible with $J$, in the sense that the associated Gram matrix satisfies $J H^{-1} J = H$. The associated Cartan involution $\theta_H$ acts on $\frg$ by $\theta_H(T) = - H^{-1} {}^t T H$.
\item If $T \in V^s(\bR)$, then there is a unique inner product $H = H_T$ compatible with $J$ such that $T$ is normal with respect to $H$ (i.e. such that $T$ commutes with its $H$-adjoint $H {}^t T H^{-1}$). We have $\cR(T) = \theta_{H_T}$.
\end{enumerate}
\end{proposition}
\begin{proof}
Let $Y_J$ denote the set of positive definite symmetric matrices $H$ satisfying $J H^{-1} J = H$. (Observe that any such matrix satisfies $\det(H) = 1$.) If $H \in Y_J$, then $\theta(X) = - H^{-1} {}^t X H$ is a Cartan involution of $\frg$. (It preserves $\frg$ because $H$, $J$ are compatible, and $B_\theta$ is positive definite because $H$ is an inner product.) Thus there is a map $Y_J \to X_G$ which is easily checked to be $G(\bbR)$-equivariant, where $x \in G(\bbR)$ acts on $H$ by $H \mapsto {}^t x^{-1}  H x^{-1} $. We need to check that this map is bijective. To see that it is surjective, note that $Y_J$ is non-empty, as it contains the identity matrix $H_0$ (otherwise said, the standard inner product on $\bbR^{2g+1}$), and $G(\bbR)$ acts transitively on $X_G$. To see that is injective, it suffices (by $G(\bbR)$-equivariance) to show that if $\theta_H = \theta_{H_0}$, then $H = H_0$. However, if $\theta_H = \theta_{H_0}$ then $- H^{-1} {}^t X H = - {}^t X$ for all $X \in \frg$, so $H$ is scalar (by Schur's lemma). Since $H$ is positive definite and $\det(H) = 1$, this forces $H = H_0$.

Since the Lie bracket on $\frs \frl_{2g+1}$ is given by $[X, Y] = XY - YX$, the second part of the proposition is asking for a unique $H$ such that $[T, \theta_H(T)] = 0$. The existence and uniqueness therefore follows from Theorem \ref{thm_special_cartan_involution}.
\end{proof}
To complete Proposition \ref{prop_explicit_reduction_theory}, we explain how to compute $H$ explicitly in terms of $T$. Supposing that $H$ exists, we extend it to a Hermitian inner product on $\bbC^{2g+1}$. Since the characteristic polynomial of $T$ has no repeated roots, $T$ is diagonalisable, and we can find $P \in \GL_{2g+1}(\bbC)$ such that $P^{-1} T P$ is diagonal. Since $T$ is normal with respect to $H$, its eigenvectors are orthogonal and so $D = {}^t P H \overline{P}$ is also diagonal, with positive real diagonal entries. Since $T$ is also self-adjoint with respect to $J$, ${}^t P J P$ is also diagonal, and thus ${}^t P H \overline{P}$, ${}^t P J P$ commute. Starting with the identity $H = J H^{-1} J $, we get
\[ {}^t P H \overline{P} = {}^t P J H^{-1} J \overline{P} = {}^t P J P P^{-1} H^{-1} {}^t \overline{P}^{-1} {}^t \overline{P} J \overline{P} = {}^t P J P ({}^t \overline{P} H P)^{-1} {}^t\overline{P} J \overline{P}, \]
hence 
\[ D^2 =  {}^t P H \overline{P}{}^t \overline{P} H P = {}^t P J P {}^t\overline{P} J \overline{P}. \]
This characterizes $H$ uniquely as 
\begin{equation}\label{eqn_reduction_covariant} H = {}^t P^{-1} ( {}^t P J P {}^t\overline{P} J \overline{P} )^{1/2} \overline{P}^{-1}. 
\end{equation}
where $P \in \GL_{2g+1}(\bC)$ is such that $P^{-1} T P$ is diagonal.

\section{A reduction algorithm for $\SO(J)$}\label{sec_reduction_algorithm}

In this section, we will describe a reduction algorithm for the symmetric space $X_G$ associated to the group $G = \SO(J)$ introduced in \S \ref{sec_example_grading}. All of the action takes place on the group $G$: the ambient grading plays no role. To orient the reader, we note that we will have to consider three symmetric bilinear forms on $\bbR^{2g+1}$ simultaneously. First, the form $J$ defining the group $G$; second, the standard inner product $H_0$ on $\bbR^{2g+1}$, which defines a base point in $X_G$; and third, a second inner product $H$, compatible with $J$, which we hope to bring closer to $H_0$ using the action of the group $G(\bbZ)$.

 We continue to define
\[ G = \mathrm{SO}(J) = \{ x \in \SL_{2g+1} \mid {}^t x J x = J \}, \]
which we think of as a group scheme over $\bbZ$. We write $A, N \subset G$ for the subgroups consisting of the diagonal and unipotent upper-triangular matrices, respectively. Then $G_{\bbQ}$ is reductive, and $A_\bbQ$, $A_\bbQ N_\bbQ$ are a maximal torus and Borel subgroup. 

Let $H_0$ denote the inner product on $\bbR^{2g+1}$ with respect to which $e_{-g}, \dots, e_g$ is an orthonormal basis. We define $K \subset G(\bbR)$ to be the subgroup of matrices that preserve $H_0$.
\begin{proposition}
\begin{enumerate}
\item $H_0$ is compatible with $J$, and $K$ is a maximal compact subgroup of $G(\bbR)$.
\item The product map $K \times A(\bbR)^0 \times N(\bbR) \to G(\bbR)$, $(k, a, n) \mapsto kan$, is a diffeomorphism.
\item If $H$ is an inner product which is compatible with $J$, then we can find a unique pair $(a, n) \in A(\bbR)^0 \times N(\bbR)$ such that $H = {}^t(an) H_0 an$. 
\end{enumerate}
\end{proposition}
\begin{proof}
Looking at Proposition \ref{prop_explicit_reduction_theory} we see that $H_0$ is compatible with $J$ and $\theta_0(x) = H_0^{-1} {}^t x^{-1} H_0 = {}^t x^{-1}$ is a Cartan involution of $G_\bR$. Therefore $K = G(\bbR)^{\theta_0}$ is a maximal compact subgroup of $G(\bbR)$. The second part of the proposition is a statement of the Iwasawa decomposition of $G(\bbR)$, determined by the data of $\theta_0$, $A$, and $N$ \cite[Theorem 7.31]{Kna02}. The third part of the proposition follows from the uniqueness of the Iwasawa decomposition and the fact that $G(\bbR)$ acts transitively on $X_G$ with $\mathrm{Stab}_{G(\bbR)}(\theta_0) = K$.
\end{proof}
The components of the  Iwasawa decomposition can be computed using the Gram--Schmidt process. Let $e^\ast_{-g}, \dots, e^\ast_g$ denote the result of carrying out the Gram--Schmidt orthogonalization process on the basis $e_{-g}, \dots, e_g$ with respect to the inner product defined by $H$. Thus we have formulae
\[ e_{-g} = e^\ast_{-g},  \dots, e_{j} = e_{j}^\ast + \sum_{i = -g}^{j-1} \mu_{i, j} e_{i}^\ast, \dots, e_g = e_g^\ast + \sum_{i = -g}^{g-1} \mu_{i, g} e_i^\ast. \]
with \[ \mu_{i, j} = ( e_j, e_i^\ast )_H / (e_i^\ast, e_i^\ast)_H. \]
In particular, we take $\mu_{j, j} = 1$ and $\mu_{i, j} = 0$ if $i > j$.
\begin{lemma}
Let $n = n_H = (\mu_{i, j})_{-g \leq i, j \leq g}$ and $a = a_H = \diag( \| e^\ast_i\|_{H})_{-g \leq i \leq g}$. Then $a \in A(\bbR)^0$, $n \in N(\bbR)$, and $H = {}^t (an) an$. 
\end{lemma} 
\begin{proof}
Let $\theta_1$ denote the Cartan involution $\theta_1 : \SL_{2g+1} \to \SL_{2g+1}$, $g \mapsto {}^t g^{-1}$. Let $K_1 = \SL_{2g+1}(\bbR)^{\theta_1}$ and let $A_1, N_1$ be the subgroups of diagonal and unipotent upper-triangular matrices in $\SL_{2g+1}$, respectively. The Iwasawa decomposition for $\SL_{2g+1}$ is the statement that the product map 
\[ K_1 \times A_1(\bbR)^0 \times N_1(\bbR) \to \SL_{2g+1}(\bbR) \]
 is a diffeomorphism. It is a standard fact that the Gram--Schmidt process gives the Iwasawa decomposition for $\SL_{2g+1}$, in the sense that $H = {}^t (an) an$ and that these are the unique elements $a \in A_1(\bbR)^0$, $n \in N_1(\bbR)$ with this property. We need to explain why the assumption that $H$, $J$ are compatible implies that in fact $a$, $n$ lie in $G(\bbR)$. However, we have $K \leq K_1$, $A \leq A_1$, and $N \leq N_1$, by construction, so the existence of the Iwasawa decomposition for $G$ and the uniqueness for $\SL_{2g+1}$ implies that the Gram--Schmidt process for $\SL_{2g+1}$ must be compatible with the Iwasawa decomposition for $G$. 
\end{proof}

\begin{lemma}\label{lem_size_reduced}
For any $n= (n_{ij}) \in N(\bbR)$, there exists $m \in N(\bbZ)$ satisfying the following conditions:
\begin{enumerate}
\item $|(nm)_{ij}| \leq 1/2$ for all $-g \leq i < j$, $j = -g, \dots, -1$.
\item $|(nm)_{i, 0}| \leq 1$ for all $-g \leq i \leq -1$.
\item $|(nm)_{i, j}| \leq 1/2$ for all $-g \leq i < -j$, $j = 1, \dots, g$.
\end{enumerate}
\end{lemma}
\begin{proof}
Right multiplication by $m$ corresponds to performing column operations on $n$. We can check that the following column operations on $n$ are induced by elements of $N(\bbZ)$:
\begin{itemize}
\item $A(i, j, q)$: For $q \in \bbZ$, $-g \leq j \leq -1$, and $-g \leq i < j$, add $q$ times column $i$ to column $j$ and $-q$ times column $-j$ to column $-i$.
\item $B(i, q)$: $q \in 2 \bbZ$, $-g \leq i \leq -1$, add $q$ times column $i$ to column $0$, and add $-q$ times column $0$ and $-q^2/2$ times column $i$ to column $-i$.
\item $C(i, j, q)$: $q \in \bbZ$, $1 \leq j \leq g$, and $-g \leq i \leq -j$, add $q$ times column $i$ to column $j$ and $-q$ times column $-j$ to column $-i$. 
\end{itemize}
We therefore carry out column operations as follows in order to satisfy the conditions of the lemma (noting that the order of operations is chosen so that once we have forced a given $n_{ij}$ to satisfy the conditions of the lemma, its value will not be changed by later operations):
\begin{enumerate}
\item For each $j = -g, \dots -1$, then for each $i = j-1, j-2, \dots, -g$, let $q$ denote the closest integer to $n_{ij}$ and do $A(i, j, -q)$.
\item For each $i = -1, \dots, -g$, let $q$ denote the closest even integer to $n_{i0}$ and do $B(i, -q)$.
\item For each $j = 1, \dots, g$, then for each $i = -j-1, \dots, -g$, let $q$ denote the closest integer to $n_{ij}$ and do $C(i, j, -q)$. 
\end{enumerate}
\end{proof}
The following definition is the analogue in our context of \cite[(1.4), (1.5)]{Len82}. 
\begin{defn}
Let $\delta \in (1/2, 1)$, and let $H$ be an inner product on $\bbR^{2g+1}$ compatible with $J$. We say that $H$ is $\delta$-reduced if the following conditions are satisfied:
\begin{enumerate}
\item $n_H$ satisfies the conditions of Lemma \ref{lem_size_reduced}.
\item For each $i = -g+1, \dots, -1$, we have $\| e_i^\ast + \mu_{i-1, i} e_{i-1}^\ast \|_H^2 \geq \delta \| e_{i-1}^\ast \|_H^2$.
\item $\| e_1^\ast + \mu_{0,1} e_0^\ast + \mu_{-1, 1} e_{-1}^\ast \|_H^2 \geq \delta^2 \| e_{-1}^\ast \|_H^2$.
\end{enumerate}
\end{defn}
\begin{remark}
The above conditions may be reformulated in terms of the matrices $n_H$, $a_H$ as follows (writing $a_H = \diag (a_{-g}, \dots, a_g)$):
\begin{enumerate} 
\item $n_H$ satisfies the conditions of Lemma \ref{lem_size_reduced}.
\item For each $i = -g +1, \dots, -1$, we have $a_{i}^2 + \mu_{i-1, i}^2 a_{i-1}^2 \geq \delta a_{i-1}^2$ (as the $e_j^\ast$ are pairwise orthogonal).
\item $a_{-1}^{-2} + \frac{1}{2}\mu_{-1, 0}^2  \geq \delta$ (noting that as $n_H, a_H \in G(\bbR)$, we have the relations $a_{-1} a_1 = a_0 = 1$, $\mu_{0, 1} = -\mu_{-1, 0}$, and $\mu_{-1, 1} = - \mu_{-1, 0}^2 / 2$). 
\end{enumerate} 
\end{remark}
We note that these conditions involve only about half the entries of the matrices $n_H, a_H$. However, because  these matrices lie in the group $\SO_{2g+1}(\bR)$, these entries determine the remaining ones. 

The following lemma shows that if $H$ is $\delta$-reduced, then it lies in a Siegel set. We recall (see \cite[Theorem 15.4]{Bor19}) that Siegel sets are fundamental sets for the action of the arithmetic group $G(\bbZ)$, in the sense defined in the introduction to this paper. This lemma can be used to show that if $H$ is $\delta$-reduced, then it has other desirable properties. For example, arguing as in the proof of \cite[Proposition 1.6]{Len82}, one can show that there is a constant $c = c_{g, \delta} > 0$ such that if $H$ is $\delta$-reduced, then $\prod_{j = -g}^g \| e_j \|^2_H \leq c$. (We thank the referee for this remark.)
\begin{lemma}
If $H$ is $\delta$-reduced, then $a_H n_H$ lies in the Siegel set $\mathfrak{S}_\delta = K A_\delta N_c$ defined as follows:
\begin{multline*} A_\delta = \{ a = (a_{-g}, \dots, a_g) \in A(\bbR)^0 \mid \\ \forall i = -g, \dots, -2, a_i / a_{i+1} \leq (\delta - 1/4)^{-1/2}; \\a_{-1} \leq (\delta - 1/2)^{-1/2} \},
\end{multline*} 
\begin{multline*} N_c = \{ n \in N(\bbR) \mid \\ \forall j = -g, \dots, -1, -g \leq i < j, \max(| n_{ij} |, |n_{i,-j}| )\leq 1/2; \\ \forall i = -g, \dots, -1, |n_{i0}| \leq 1 \}.
\end{multline*}
\end{lemma}
\begin{proof}
Since $n_H$ lies in $N_c$ by definition of being $\delta$-reduced, we just need to check that $a_H$ satisfies the required inequalities. Re-arranging gives
\[ a_{i+1}^2 / a_{i}^2 \geq \delta - \mu_{i, i+1}^2 \geq \delta - 1/4\]
for each $i = -g, \dots, -2$, and
\[ a_{1}^2 \geq \delta - \mu_{-1, 0}^2 / 2 \geq \delta - 1/2, \]
as required. 
\end{proof}
The Weyl group $W(G, A)$ is generated by its simple reflections (with respect to the choice of positive roots determined by $N$). These simple reflections can be represented by the following elements $s_{-g}, \dots, s_{-1}$ of $G(\bbZ)$, which we define by their action on basis vectors:
\begin{itemize}
\item For each $i = -g, \dots, -2$, $s_i$ sends $e_i$ to $e_{i+1}$, $e_{i+1}$ to $e_i$, $e_{-i}$ to $e_{-(i+1)}$, $e_{-(i+1)}$ to $e_{-i}$, and fixes the remaining basis vectors.
\item $s_{-1}$ sends $e_{-1}$ to $e_1$, $e_1$ to $e_{-1}$, $e_0$ to $-e_0$, and fixes the remaining basis vectors.
\end{itemize}
Here finally is an algorithm which, starting with any $H$ as above, finds $\gamma \in G(\bbZ)$ such that ${}^t \gamma H \gamma$ is $\delta$-reduced:
\begin{enumerate}
\item Let $\gamma_1 = 1$, $H_1 = H$.
\item Set $H_1 := {}^t \gamma_1 H \gamma_1$. Compute the matrices $n_{H_1} = (\mu_{i, j}), a_{H_1} = \diag(a_{-g}, \dots, a_g)$ using Gram--Schmidt.  
\item If $n_{H_1}$ satisfies the conditions of Lemma \ref{lem_size_reduced}, then proceed to Step 4. Else, find $\gamma_n$ such that $n_{H_1} \gamma_n$ satisfies the conditions of Lemma \ref{lem_size_reduced}, set $\gamma_1 := \gamma_1 \gamma_n$, and go to Step 2. 
\item For each $i = -g, \dots, -2$: if $a_{i+1}^2 + \mu_{i, i+1}^2 a_{i}^2 < \delta a_{i}^2$, set $\gamma_1 := \gamma_1  s_i$ and go to Step 2. 
\item If $a_{-1}^{-2} + \frac{1}{2} \mu_{-1, 0}^2 < \delta$, set $\gamma_1 := \gamma_1 s_{-1}$ and go to Step 2.
\item Return $\gamma_1$. 
\end{enumerate}
In practice, it is more efficient to run Gram--Schmidt once (for the matrix $H$) and then to update the matrices $n_{H_1}$, $a_{H_1}$ at each step. Explicit formulae describing how to do this are given in \cite{Rom24}. 
\begin{proposition}
The above algorithm always terminates. 
\end{proposition}
\begin{proof}
It suffices to show that the simple reflections $s_i$ are applied only finitely many times. At each instance of Step 2 we are given a basis $f_i = \gamma_1 e_i$ of $\bbZ^{2g+1}$, to which we associate the flag $F_i = \oplus_{k=-g}^i \bbZ f_k$, $i = -g, \dots, g$. If $i = -g, \dots, -2$, then acting by $s_i$ changes $F_i$ and $F_{-(i+1)}$ but leaves the remaining $F_j$ unchanged. Acting by $s_{-1}$ changes $F_{-1}$ and $F_0$ and leaves the remaining $F_j$ unchanged. In particular, among $F_{-g}, \dots, F_{-1}$, $s_i$ changes only $F_i$.

The inner product $H$ on $\bbR^{2g+1}$ determines one on $\wedge^i \bbR^{2g+1}$ for each $i = 1, \dots, 2g+1$ by the usual formula $( x_1 \wedge \dots \wedge x_i, y_1 \wedge \dots \wedge y_i )_H = \det( (x_j, y_k)_H )$. Since $\wedge^i \bbZ^{2g+1}$ contains finitely many vectors of bounded norm, it suffices to show that applying the simple reflection $s_i$ decreases the squared norm of the vector (which, up to sign, depends only on $F_i$) $f_{-g} \wedge \dots \wedge f_{i} \in \wedge^{i+g+1} \bbZ^{2g+1}$. A calculation shows that this squared norm is multiplied by a scalar which is strictly less than $\delta$. As $\delta < 1$, this completes the proof.
\end{proof}
\section{Example: Reducing a self-adjoint linear operator}\label{sec_example_computation}

We now show how to combine the theory of the previous two sections in an example. Take $g = 3$, and consider the group $G = \SO(J) \subset \SL_7$ and representation $V \subset \frs\frl_7$ associated to the $\bbZ / 2 \bbZ$-grading of $\frs\frl_7$ defined in \S \ref{sec_example_grading}. We consider the genus 3 hyperelliptic curve 
\[ C_f : y^2 = f(x) =    x^7- x^5 - 2 x^4 - x^3+ 5 x -5 \]
This curve has the integral point $P = (14, 10237)$. Associated to $P$ is the matrix
\[ T = \left(
\begin{array}{ccccccc}
 -14 & 1 & 0 & 0 & 0 & 0 & 0 \\
 -195 & 0 & 1 & 0 & 0 & 0 & 0 \\
 -2728 & 0 & 7 & 0 & -1 & 0 & 0 \\
 -10237 & 0 & 0 & 14 & 0 & 0 & 0 \\
 19095 & -6 & -48 & 0 & 7 & 1 & 0 \\
 1546 & -26 & -6 & 0 & 0 & 0 & 1 \\
 390 & 1546 & 19095 & -10237 & -2728 & -195 & -14 \\
\end{array}
\right) \in V(\bZ) \]
of characteristic polynomial $f(x)$, which we obtain e.g. from the subregular Slodowy slice studied in \cite{Tho13}. We compute the reduction covariant $H = \cR(T) \in X_{\SO(J)}$ using the formula (\ref{eqn_reduction_covariant}). We obtain numerically
\[ \Scale[0.8]{ H =  \left(
\begin{array}{ccccccc}
 3.74708 & 53.7691 & 750.242 & 2813.43 & -5244.78 & -421.526 & -47.2448 \\
 53.7691 & 776.143 & 10830.1 & 40612.6 & -75708.6 & -6080.03 & -681.676 \\
 750.242 & 10830.1 & 151130. & 566729. & -1.05648\times 10^6 & -84842.6 & -9520.71 \\
 2813.43 & 40612.6 & 566729. & 2.12521\times 10^6 & -3.96175\times 10^6 & -318157. & -35704.6 \\
 -5244.78 & -75708.6 & -1.05648\times 10^6 & -3.96175\times 10^6 & 7.38537\times 10^6 & 593097. & 66564.2 \\
 -421.526 & -6080.03 & -84842.6 & -318157. & 593097. & 47660.8 & 5338.34 \\
 -47.2448 & -681.676 & -9520.71 & -35704.6 & 66564.2 & 5338.34 & 660.273 \\
\end{array}
\right).} \]
We then apply the algorithm of \S \ref{sec_reduction_algorithm} to $H$ with $\delta = 0.9$ to obtain an element $\gamma_1 \in G(\bbZ)$ such that ${}^t \gamma_1 H \gamma_1$ is $\delta$-reduced. The resulting matrix is 
\[ \gamma_1 = \left(
\begin{array}{ccccccc}
 -2 & -1 & -15 & 2 & -3 & -5 & 26 \\
 0 & -8 & -117 & 12 & 0 & 9 & 203 \\
 -8 & -104 & -1462 & 180 & -16 & 184 & 2557 \\
 4 & 56 & 784 & -97 & 8 & -98 & -1372 \\
 1 & 15 & 209 & -26 & 2 & -26 & -366 \\
 0 & 1 & 15 & -2 & 0 & -2 & -26 \\
 0 & 0 & 1 & 0 & 0 & 0 & -2 \\
\end{array}
\right), \]
yielding
\[ \gamma_1^{-1} T \gamma_1 = \left(
\begin{array}{ccccccc}
 0 & 0 & -1 & 2 & 2 & -2 & 3 \\
 1 & 0 & 1 & 0 & 0 & 0 & -2 \\
 0 & 1 & 1 & -2 & 0 & 0 & 2 \\
 0 & 0 & 1 & -2 & -2 & 0 & 2 \\
 0 & 0 & 0 & 1 & 1 & 1 & -1 \\
 0 & 0 & 0 & 0 & 1 & 0 & 0 \\
 0 & 0 & 0 & 0 & 0 & 1 & 0 \\
\end{array}
\right). \]
\begin{remark} 
Since the lower left-hand $3 \times 3$ submatrix of $\gamma_1^{-1} T \gamma_1$ is 0, we see that the orbit of $T$ is distinguished, in the sense of \cite{Bha13}. It follows that the divisor class of $P - P_\infty$ is divisible by 2 in $J_f(\bQ)$ (where $P_\infty$ is the point at infinity of $C_f$ and $J_f$ is the Jacobian of $C_f$). This could also be checked using the usual theory of 2-descent (see e.g. \cite{Sch95}), by noting that the class of $14-x$ in the group $L_f^\times /  (L_f^\times)^2$ (where $L_f = \bbQ[x] / (f(x))$) is trivial. 
\end{remark} 

\bibliographystyle{alpha}
\bibliography{reduction}

\begin{thebibliography}{RLYG12}

\bibitem[AT18]{Ada18}
Jeffrey Adams and Olivier Ta\"{\i}bi.
\newblock Galois and {C}artan cohomology of real groups.
\newblock {\em Duke Math. J.}, 167(6):1057--1097, 2018.

\bibitem[BG13]{Bha13}
Manjul Bhargava and Benedict~H. Gross.
\newblock The average size of the 2-{S}elmer group of {J}acobians of
  hyperelliptic curves having a rational {W}eierstrass point.
\newblock In {\em Automorphic representations and {$L$}-functions}, volume~22
  of {\em Tata Inst. Fundam. Res. Stud. Math.}, pages 23--91. Tata Inst. Fund.
  Res., Mumbai, 2013.

\bibitem[BHC62]{Bor62}
Armand Borel and Harish-Chandra.
\newblock Arithmetic subgroups of algebraic groups.
\newblock {\em Ann. of Math. (2)}, 75:485--535, 1962.

\bibitem[Bor19]{Bor19}
Armand Borel.
\newblock {\em Introduction to arithmetic groups}, volume~73 of {\em University
  Lecture Series}.
\newblock American Mathematical Society, Providence, RI, 2019.
\newblock Translated from the 1969 French original by Lam Laurent Pham, Edited
  and with a preface by Dave Witte Morris.

\bibitem[CFS10]{Cre10}
John~E. Cremona, Tom~A. Fisher, and Michael Stoll.
\newblock Minimisation and reduction of 2-, 3- and 4-coverings of elliptic
  curves.
\newblock {\em Algebra Number Theory}, 4(6):763--820, 2010.

\bibitem[Fis13]{Fis13}
Tom Fisher.
\newblock Minimisation and reduction of 5-coverings of elliptic curves.
\newblock {\em Algebra Number Theory}, 7(5):1179--1205, 2013.

\bibitem[Kna02]{Kna02}
Anthony~W. Knapp.
\newblock {\em Lie groups beyond an introduction}, volume 140 of {\em Progress
  in Mathematics}.
\newblock Birkh\"{a}user Boston, Inc., Boston, MA, second edition, 2002.

\bibitem[Lag22]{Lag22}
Jef Laga.
\newblock The average size of the 2-{S}elmer group of a family of
  non-hyperelliptic curves of genus 3.
\newblock {\em Algebra Number Theory}, 16(5):1161--1212, 2022.

\bibitem[LLL82]{Len82}
A.~K. Lenstra, H.~W. Lenstra, Jr., and L.~Lov\'{a}sz.
\newblock Factoring polynomials with rational coefficients.
\newblock {\em Math. Ann.}, 261(4):515--534, 1982.

\bibitem[Mos55]{Mos55}
G.~D. Mostow.
\newblock Self-adjoint groups.
\newblock {\em Ann. of Math. (2)}, 62:44--55, 1955.

\bibitem[RLYG12]{Ree12}
Mark Reeder, Paul Levy, Jiu-Kang Yu, and Benedict~H. Gross.
\newblock Gradings of positive rank on simple {L}ie algebras.
\newblock {\em Transform. Groups}, 17(4):1123--1190, 2012.

\bibitem[RT24]{Rom24}
Beth Romano and Jack~A. Thorne.
\newblock An {LLL} algorithm with symmetries.
\newblock Preprint, 2024.

\bibitem[SC03]{Sto03}
Michael Stoll and John~E. Cremona.
\newblock On the reduction theory of binary forms.
\newblock {\em J. Reine Angew. Math.}, 565:79--99, 2003.

\bibitem[Sch95]{Sch95}
Edward~F. Schaefer.
\newblock {$2$}-descent on the {J}acobians of hyperelliptic curves.
\newblock {\em J. Number Theory}, 51(2):219--232, 1995.

\bibitem[Sha19]{Sha19}
Ananth~N. Shankar.
\newblock 2-{S}elmer groups of hyperelliptic curves with marked points.
\newblock {\em Trans. Amer. Math. Soc.}, 372(1):267--304, 2019.

\bibitem[Tho13]{Tho13}
Jack~A. Thorne.
\newblock Vinberg's representations and arithmetic invariant theory.
\newblock {\em Algebra Number Theory}, 7(9):2331--2368, 2013.

\bibitem[Tho15]{Tho15}
Jack~A. Thorne.
\newblock {$E_6$} and the arithmetic of a family of non-hyperelliptic curves of
  genus 3.
\newblock {\em Forum Math. Pi}, 3:e1, 41, 2015.

\bibitem[Vin76]{Vin76}
\`E.~B. Vinberg.
\newblock The {W}eyl group of a graded {L}ie algebra.
\newblock {\em Izv. Akad. Nauk SSSR Ser. Mat.}, 40(3):488--526, 709, 1976.

\end{thebibliography}

\end{document}